\documentclass[12pt]{article}
\usepackage{amsmath}
\usepackage{amsfonts}
\usepackage{mathrsfs}
\usepackage{amssymb}
\usepackage[usenames]{color}
\usepackage[numbers,sort&compress]{natbib}
\usepackage{subfigure}
\usepackage{bm}
\usepackage[all]{xy}
\usepackage{setspace}
\usepackage{multirow}

\def\cl{\centerline}
\def\ni{\noindent}
\def\vs{\vspace*}

\def\a{\alpha}

\def\D{\Delta}

\def\CVir{{\mathfrak {Vir}}}

\def\Z{\mathbb{Z}}

\def\C{\mathbb{C}}

\def\QED{\hfill$\Box$}

\numberwithin{equation}{section}
\newtheorem{theo}{Theorem}[section]
\newtheorem{defi}[theo]{Definition}

\newtheorem{lemm}[theo]{Lemma}
\newtheorem{rema}[theo]{Remark}
\newtheorem{prop}[theo]{Proposition}

\newtheorem{exam}{Example}

\begin{document}
\begin{center}
{\bf\Large Virasoro elements in conformal algebras}
\footnote{
$^{*}$Corresponding author: chgxia@cumt.edu.cn (C.~Xia).}
\end{center}
\vspace{15pt}

\cl{Yucai Su$^{\,\dag}$, \ Chunguang Xia$^{\,*,\,\ddag,\,\S}$}
\vspace{20pt}

\cl{\small $^{\,\dag}$Department of Mathematics, Jimei University, Xiamen 361021, China}
\vspace{8pt}

\cl{\small $^{\,\ddag}$School of Mathematics, China University of Mining and Technology, Xuzhou 221116, China}
\vspace{8pt}

\cl{\small $^{\,\S}$Jiangsu Center for Applied Mathematics (CUMT), Xuzhou 221116, China}
\vspace{8pt}

\cl{\small E-mails: ycsu@tongji.edu.cn, chgxia@cumt.edu.cn}
\vspace{8pt}

\footnotesize
\noindent{{\bf Abstract.}
It is well known that Virasoro elements play an extremely important role in the theory of vertex operator algebras, 
and likewise a very significant role in the field of conformal algebras.
In this paper, we introduce the notion of a regular action in the category of conformal modules over Lie conformal algebras with Virasoro elements.
We show that a finite conformal module over the general conformal algebra $\mathfrak{gc}_1$
(resp., $\mathfrak{gc}_N$ with $N\ge2$) is semisimple if and only if there exists a pair of different Virasoro elements
(resp., canonical Virasoro elements) with regular actions.
 Along the way to finding semisimplicity criteria, we also discuss the classification of Virasoro elements of $\mathfrak{gc}_N$ in-depth, leading us
to construct a huge number of new Virasoro conformal modules.
\vs{5pt}

\ni{\bf Key words:} Virasoro element; general conformal algebra; conformal module; semisimplicity criteria

\ni{\it Mathematics Subject Classification (2020):} 17B10; 17B65; 17B68; 17B69.}

\small
\begin{spacing}{1}
\renewcommand{\baselinestretch}{1}
\tableofcontents
\end{spacing}

\small
\section{Introduction}
The investigation of Lie conformal (super)algebras and their representation theory was initiated by Kac \cite{K1996-98},
motivated by their deep relationships with the $2$-dimensional quantum ﬁeld theory, vertex algebras, and infinite-dimensional Lie (super)algebras.

The simplest but fundamental finite simple Lie conformal algebra is the Virasoro conformal algebra $\CVir$.
The most important infinite simple Lie conformal algebra is the general conformal algebra $\mathfrak{gc}_N$, which is the conformal analogue of
the general Lie algebra $\mathfrak{gl}_N$.
These two simple Lie conformal algebras play the central role in the whole Lie conformal (super)algebras theory.
Irreducible finite conformal modules over $\CVir$ and $\mathfrak{gc}_N$ were respectively classified in \cite{CK1997} and
a work due to Kac, Radul and Wakimoto (see also \cite{K1999,Re2006,BM2013} for details) around in 1997.

However, for general finite conformal modules over $\CVir$ and $\mathfrak{gc}_N$,
there do not exist conformal analogues of Weyl's complete reducibility theorem \cite{K1996-98,CKW1998,BKV1999,S2004,Re2006},
although $\CVir$ and $\mathfrak{gc}_N$ are simple.
Hence, finding precise and verifiable semisimplicity criteria for these conformal modules is a problem of great interest \cite{K1999},
which has remained open for nearly thirty years.

In this paper, by abstracting a common feature of irreducible conformal modules over $\CVir$, $\mathfrak{gc}_N$ and some other typical
Lie conformal algebras, we propose the notion of a {\it regular action} in the category of conformal modules over Lie conformal algebras 
with Virasoro elements (Definition~\ref{def-ra}). Under this definition, it is clear that a finite conformal module
over $\CVir$ is semisimple if and only if the action of the unique Virasoro element of $\CVir$ is regular and all conformal weights, 
with respect to some $\C[\partial]$-basis, are non-zero (Theorem~\ref{thm-vir}).
For $\mathfrak{gc}_N$, we introduce the notion of a {\it canonical Virasoro element},
which is defined to be the homomorphism image of a Virasoro element of $\mathfrak{gc}_1$ under the
canonical embedding map from $\mathfrak{gc}_1$ into $\mathfrak{gc}_N$.
We prove the following semisimplicity criteria.


\vs{6pt}
\noindent {\bf Main Theorem.}
\label{main-thm}
{\it Let $V$ be a finite conformal module over $\mathfrak{gc}_N$.
\begin{itemize}\parskip-4pt
\item[\rm(1)] If $N=1$, then $V$ is semisimple if and only if there exists a pair of different Virasoro elements,
whose actions on a $\C[\partial]$-basis of $V$ are regular.
\item[\rm(2)] If $N\ge2$, then $V$ is semisimple if and only if there exists a pair of different canonical Virasoro elements,
whose actions on a $\C[\partial]$-basis of $V$ are regular.
\end{itemize}}

The key step in the proof of part (1) of the Main Theorem is the observation that $\mathfrak{gc}_1$ contains a Heisenberg-Virasoro conformal subalgebra,
whose representation theory 
will be used.
While in the proof of part (2) of the Main Theorem, the rigidities of (anti-)homomorphisms between matrix algebras,
guaranteed by the Skolem-Noether theorem and the semisimplicity of representations of matrix algebras, play a crucial role.

Along the way to finding semisimplicity criterion, we are also concerned with another related interesting problem:
classify Virasoro elements of $\mathfrak{gc}_N$. For $N=1$, the classification was stated earlier in \cite{Re2006}
(see also \cite{DeK2002,BKL2003-2} and Theorem~\ref{prop-vir-ele-1}); here we give a complete proof.
For $N\ge2$, we introduce the notion of a {\it standard Virasoro element} (Definition~\ref{def-sv}),
give an explicit classification of those of degree one (Theorem~\ref{prop-vir-ele-N}), and construct those of
higher degree and non-standard Virasoro elements (Propositions~\ref{std-vir-ele-higher} and \ref{non-std-vir-ele}),
leading us to construct a huge number of new Virasoro conformal modules (Remark~\ref{new-vir-modules}).
Note the fact that a description of all Virasoro elements in $\mathfrak{gc}_2$ follows from \cite{CKW1998} (see also \cite{DeK2002}), and that the analogous description of those 
inside $\mathfrak{gc}_3$, say, would follow from a classification of all double extensions, which is an attainable and a finite result.

This work is organized as follows.
In Section~2, we recall some basic definitions and propose the notion of a regular action.
In Section 3, we introduce the notions of canonical and standard Virasoro elements of $\mathfrak{gc}_N$.
We provide an in-depth discussion of the classification of Virasoro elements of $\mathfrak{gc}_N$ as summarized in the above paragraph.
Sections~4 and 5 are devoted to proving the Main Theorem.
Finally, in Section~6, we propose some open problems arising from our study.

\section{Preliminaries}

Throughout this paper, we work over the complex number field $\C$.

\subsection{Basic definitions}

First, we recall some basic definitions on Lie conformal algebras, see \cite{K1996-98,CK1997,CCK2021,CCK2024} for more details.

\begin{defi}\label{def-conforma-algebra}\rm
A {\it Lie conformal algebra} $R$ is a $\C[\partial]$-module endowed with a
$\C$-linear map \linebreak $R\otimes R\rightarrow \C[\lambda]\otimes R$, $a\otimes b\rightarrow [a\,{}_\lambda \,b]$,
called the {\it $\lambda$-bracket}, satisfying the axioms ($a,\,b,\,c\in R$):
\begin{equation*}
\aligned
\mbox{(C1)}&~~~~[\partial a\,{}_\lambda \,b]=-\lambda[a\,{}_\lambda\, b],\ \ \ \
[a\,{}_\lambda \,\partial b]=(\partial+\lambda)[a\,{}_\lambda\, b],\\
\mbox{(C2)}&~~~~[a\, {}_\lambda\, b]=-[b\,{}_{-\lambda-\partial}\,a],\\
\mbox{(C3)}&~~~~[a\,{}_\lambda\,[b\,{}_\mu\, c]]=[[a\,{}_\lambda\, b]\,{}_{\lambda+\mu}\, c]+[b\,{}_\mu\,[a\,{}_\lambda \,c]].
\endaligned
\end{equation*}
\end{defi}

\begin{defi}\label{def-conformal-module}\rm
A {\it conformal module} $M$ over a Lie conformal algebra $R$ is a $\C[\partial]$-module endowed with  a $\C$-linear map
$R\otimes M\rightarrow \C[\lambda]\otimes M$, $a\otimes v\rightarrow a\,{}_\lambda \,v$, called the {\it $\lambda$-action},
satisfying the properties ($a,\,b\in R$, $v\in M$):
\begin{equation*}
\aligned
\mbox{(M1)}&~~~~(\partial a)\,{}_\lambda\, v=-\lambda a\,{}_\lambda\, v,\ \ \ \ \ a{}\,{}_\lambda\, (\partial v)=(\partial+\lambda)a\,{}_\lambda\, v,\\
\mbox{(M2)}&~~~~[a\,{}_\lambda\, b]\,{}_{\lambda+\mu}\, v = a\,{}_\lambda\, (b{}\,_\mu\, v)-b\,{}_\mu\,(a\,{}_\lambda\, v).
\endaligned
\end{equation*}
\end{defi}

A Lie conformal algebra $R$ or a conformal $R$-module $M$ is called {\it finite} if it is finitely generated as a $\C[\partial]$-module.
A conformal $R$-module is called {\it simple} (or {\it irreducible}) if it contains no non-trivial submodules. 
A conformal $R$-module is called {\it semisimple} (or {\it completely reducible}) if it is a direct sum of simple (or irreducible) submodules.

For any conformal $R$-module $M$, it was shown in \cite{DK1998} that $R$ acts trivially on the torsion of $M$. In view of this, we are only concerned with
modules without $\C[\partial]$-torsion.

Denote by $M^{(n)}$ the direct sum of $n$ copies of a conformal $R$-module $M$.

\begin{defi}\label{def-conformal-dual}\rm
Let $M$ be a conformal module over a Lie conformal algebra $R$.
The {\it conformal dual} (or {\it contragredient module}) of $M$ is defined by
$$
M^*=\{f_\lambda: M\rightarrow\C[\lambda]\,\,|\,\,f_\lambda(\partial v)=\lambda f_\lambda(v),\ v\in M\}
$$
with the structures of $\C[\partial]$-module and $R$-module, respectively, given by ($a\in R$, $v\in M$):
\begin{equation*}
\aligned
\mbox{(D1)}&~~~~(\partial f)\,{}_\lambda\, (v)=-\lambda f\,{}_\lambda\, (v),\\
\mbox{(D2)}&~~~~(a\,{}_\lambda\, f)\,{}_{\mu}\, (v) = -f\,{}_{\mu-\lambda}\, (a{}\,_\lambda\, v).
\endaligned
\end{equation*}
\end{defi}

Next, we recall the most important two simple Lie conformal algebras (one is the Virasoro conformal algebra $\CVir$ and
the other is the general conformal algebra $\mathfrak{gc}_N$) and their irreducible conformal modules.

\begin{exam}\label{ex-1}\rm
The Virasoro conformal algebra $\CVir=\C[\partial]L$ is the simplest but fundamental example of finite Lie conformal algebras,
which has $\lambda$-bracket $[L\, {}_\lambda \, L]=(\partial+2\lambda) L$.
Any non-trivial free conformal module of rank one over $\CVir$ \cite{DK1998} has the form
$M^{\CVir}_{\D,\a}=\C[\partial]v$ with $\lambda$-action
\begin{equation}\label{CVir-actions}
L\,{}_\lambda\, v=(\partial+\D\lambda+\a)v,
\end{equation}
where $\D,\a\in\C$. We have $(M^{\CVir}_{\D,\a})^*\cong M^{\CVir}_{1-\D,-\a}$. Furthermore, $M^{\CVir}_{\D,\a}$ is irreducible if and only if $\D\ne 0$,
and all non-trivial finite irreducible conformal modules over $\CVir$ are of this kind \cite{CK1997}.
\end{exam}

Denote by $M_{n}$ the set of all $n\times n$ matrices, $I_n$ the identity matrix of order $n$. 
For $A\in M_{n}$, denote by $A^T$ the transpose of $A$.

\begin{exam}\label{ex-3}\rm
Let $N$ be a positive integer. The general conformal algebra $\mathfrak{gc}_N=\C[\partial,x]\otimes M_N$ has a $\C[\partial]$-generating set
$\{J^n_A:=x^n\otimes A\,|\,n\in\Z_{+},\,A\in M_N\}$ with $\lambda$-brackets
\begin{equation}\label{gcN-brackets}
[J^m_A\, {}_\lambda \, J^n_B] = \sum_{s=0}^m \binom{m}{s}(\partial+\lambda)^{s} J^{m+n-s}_{AB}- \sum_{s=0}^n \binom{n}{s}(-\lambda)^{s} J^{m+n-s}_{BA}.
\end{equation}
See \cite{DK1998} for the original definition of $\mathfrak{gc}_N$.
It is the most important infinite Lie conformal algebra, which plays the same role in the theory of Lie conformal algebras as
the general Lie algebra $\mathfrak{gl}_N$ does in the theory of Lie algebras. The $\C[\partial]$-module
$\C[\partial]^N=\C[\partial]\otimes\C^N$ becomes a conformal $\mathfrak{gc}_N$-module if we define $\lambda$-actions by
\begin{equation}\label{gcN-actions}
J^n_A\,{}_\lambda\, v=(\partial+\lambda+\a)^n Av \ \ (\mbox{resp.},\ J^n_A\,{}_\lambda\, v=-(-\partial+\a)^n A^Tv), \ \ v\in \C^N,
\end{equation}
where $\a\in\C$. Denote this module by $M^{\mathfrak{gc}_N}_{\a}$ (resp., $(M^{\mathfrak{gc}_N}_{\a})^*$).
Any non-trivial finite irreducible conformal module over $\mathfrak{gc}_N$
is isomorphic to $M^{\mathfrak{gc}_N}_{\a}$ or the conformal dual $(M^{\mathfrak{gc}_N}_{\a})^*$ of $M^{\mathfrak{gc}_N}_{\a}$.
This classification is due to Kac, Radul and Wakimoto, see also \cite{K1999,Re2006,BM2013}.
\end{exam}

\subsection{Regular actions}

A non-zero element $L$ of a Lie conformal algebra is referred to as a {\it Virasoro element} if $[L\, {}_\lambda \, L]=(\partial+2\lambda) L$.
From \eqref{CVir-actions} and \eqref{gcN-actions}, we observe that the actions of Virasoro elements have a common feature:

\begin{table}[ht]\label{tab-1}
\centering\small
\subtable[A common feature]{
\begin{tabular}{c|c|l|c}
\hline{\centering} \rule{0pt}{14pt}
Example  & Algebra & Virasoro element &  Common feature \\
\hline \rule{0pt}{14pt}
1 & $\CVir$ & $L=L$   & \multirow{2}*{
$\begin{array}{ll}
L\,{}_\lambda\, v=(\partial+\D\lambda+\a)v\\
\mbox{\ for some}\ \D,\a\in\C
\end{array}$
} \\
\cline{1-3} \rule{0pt}{14pt}
2 & $\mathfrak{gc}_N$ & $L=(a\partial+b)J^0_{I_N}+J^1_{I_N}$  &  ~ \\
\hline
\end{tabular}}
\end{table}

One can also find that the actions of Virasoro elements on irreducible conformal modules over many other
Lie conformal algebras also admit the above feature, see, e.g., \cite{SXY2018,W2019}.
Motivated by this observation, we introduce the following definition:

\begin{defi}\label{def-ra}\rm
Let $R$ be a Lie conformal algebra with a Virasoro element $L$, and $V$ be a conformal $R$-module.
The action of $L$ is {\it regular} if there exists a $\C[\partial]$-basis $Y$ of $V$ such that
\begin{equation*}\label{def-rc}
L\,{}_\lambda\, y=(\partial+\D^{L}_y\lambda+\a^{L}_y)y, \quad y\in Y,
\end{equation*}
where $\D^{L}_y, \a^{L}_y\in\C$.
We refer to $\D^{L}_y$ as a {\it conformal weight}, and $\a^{L}_y$ as a {\it conformal shift}.
\end{defi}

By Definition~\ref{def-ra} and the classification results in \cite{CK1997} (see Example~\ref{ex-1}),
one can easily describe the semisimplicity of finite conformal modules over $\CVir$.

\begin{theo}\label{thm-vir}
A finite conformal module over $\CVir$ is semisimple if and only if the action of the unique Virasoro element $L$ is regular
and all conformal weights, with respect to some $\C[\partial]$-basis, are non-zero.
\end{theo}

\section{Virasoro elements}

The notion of regular actions introduced in Definition~\ref{def-ra} is defined for Virasoro elements.
Given a Lie conformal algebra $R$, when considering a conformal $R$-module with regular actions,
we are first led to ask which elements of $R$ can serve as Virasoro elements.
To completely classify Virasoro elements of $\mathfrak{gc}_N$ is an interesting problem \cite{CKW1998,DeK2002,BKL2003-2,Re2006}.
In this section, we aim to discuss this problem as in-depth as possible.

\subsection{Virasoro elements of $\mathfrak{gc}_1$}

Replacing $J^n_A$'s by $J^n$'s in Example~\ref{ex-3}, we obtain
the conformal structures of $\mathfrak{gc}_1$ and its finite irreducible conformal modules.
We shall use the notations $J^n$'s in $\mathfrak{gc}_1$ throughout the paper.
The following statement was mentioned in \cite{Re2006}. Now we provide a rigorous proof.

\begin{theo}\label{prop-vir-ele-1}
Any Virasoro element of $\mathfrak{gc}_1$ has the form $(a\partial+b)J^0+J^1$, where $a,b\in\C$.
\end{theo}

\ni{\it Proof.}\ \
Assume that $L$ is a Virasoro element of $\mathfrak{gc}_1$. By \eqref{gcN-brackets}, $L$ must have the form
$L=f_0(\partial)J^0+f_1(\partial)J^1$, where $f_0(\partial), f_1(\partial)\in\C[\partial]$.
Since $[L\, {}_\lambda \, L]=(\partial+2\lambda) L$, a direct computation shows that
\begin{eqnarray*}
 &\!\!\!\!\!\!& \left(\lambda f_0(-\lambda)f_1(\partial+\lambda)+(\partial+\lambda)f_1(-\lambda)f_0(\partial+\lambda)\right)J^0
+(\partial+2\lambda) f_1(-\lambda)f_1(\partial+\lambda)J^1\\
   &\!\!\!=\!\!\!& (\partial+2\lambda)f_0(\partial)J^0+(\partial+2\lambda)f_1(\partial)J^1,
\end{eqnarray*}
which implies that
\begin{eqnarray}
\label{equ-J0} \lambda f_0(-\lambda)f_1(\partial+\lambda)+(\partial+\lambda)f_1(-\lambda)f_0(\partial+\lambda)&\!\!\!=\!\!\!& (\partial+2\lambda)f_0(\partial),\\
\label{equ-J1}  (\partial+2\lambda) f_1(-\lambda)f_1(\partial+\lambda) &\!\!\!=\!\!\!& (\partial+2\lambda)f_1(\partial).
\end{eqnarray}
If $f_1(\partial)=0$, then \eqref{equ-J0} implies that $f_0(\partial)=0$, and thus $L=0$, a contradiction.
Hence, $f_1(\partial)\ne0$. Comparing the degrees in $\lambda$ on both sides of \eqref{equ-J1}, we see that $f_1(\partial)=1$. Then \eqref{equ-J0} becomes
$$
\lambda f_0(-\lambda)+(\partial+\lambda)f_0(\partial+\lambda) = (\partial+2\lambda)f_0(\partial).
$$
Rewrite the above equation as
$$
\partial\left(\frac{f_0(\partial+\lambda)-f_0(\partial)}{\lambda}\right)=2f_0(\partial)-f_0(-\lambda)-f_0(\partial+\lambda).
$$
Taking $\lambda\rightarrow 0$, we obtain $\partial \frac{d}{d\partial}f_0(\partial)=f_0(\partial)-f_0(0)$.
This forces $f_0(\partial)$ to be of the form $f_0(\partial)=a\partial+b$, where $a,b\in\C$
[Note: the above {\it standard  calculus techniques} are frequently used in the representation 
theory of Lie conformal (super)algebras].
Hence, $L=(a\partial+b)J^0+J^1$. Conversely, the element $(a\partial+b)J^0+J^1$ is indeed a Virasoro element of $\mathfrak{gc}_1$.
This completes the proof.
\QED

\subsection{Canonical and standard Virasoro elements of $\mathfrak{gc}_N$}
\label{sec-cs}

Consider the canonical embedding map
$$
\Pi: \mathfrak{gc}_1\rightarrow\mathfrak{gc}_N, \quad J^n\mapsto J^n_{I_N}.
$$
For $a,b\in\C$, the homomorphism image of $L_{a,b}:=(a\partial+b)J^0+J^1\in \mathfrak{gc}_1$ is $\Pi(L_{a,b})=$\linebreak $(a\partial+b)J^0_{I_N}+J^1_{I_N}\in \mathfrak{gc}_N$.
By Theorem~\ref{prop-vir-ele-1}, $\Pi(L_{a,b})$ is a Virasoro element of $\mathfrak{gc}_N$, which has been listed in Table~1.
We refer to $\Pi(L_{a,b})$ as a {\it canonical Virasoro element} of $\mathfrak{gc}_N$.

Next, we study non-obvious Virasoro elements of $\mathfrak{gc}_N$.
Denote by $E_{ij}\in M_N$ the matrix whose $(i,j)$-entry is one and all others are zero.
Then
$
\{J^s_{E_{ij}}\in\mathfrak{gc}_N\,|\,s\in\Z_+, 1\le i,j\le N\}
$
is a $\C[\partial]$-basis of $\mathfrak{gc}_N$. Any non-zero element $g\in \mathfrak{gc}_N$ can be written as
$$
g=\sum_{s=0}^{n}\sum_{i,j=1}^N f_{ij}^{[s]}(\partial)J^s_{E_{ij}},
$$
where $f_{ij}^{[s]}(\partial)\in\C[\partial]$, and there exist $i_0$ and $j_0$ such that $ f_{i_0j_0}^{[n]}(\partial)\ne0$.
We refer to $n$ as the {\it degree} of $g$, and
$$
P_g^{[s]}(\partial)=\left(
               \begin{array}{ccc}
                 f_{11}^{[s]}(\partial) & \cdots & f_{1N}^{[s]}(\partial) \\
                 \cdots & \cdots & \cdots \\
                 f_{N1}^{[s]}(\partial) & \cdots & f_{NN}^{[s]}(\partial)  \\
               \end{array}
             \right)\in M_N(\C[\partial]),\ \  0\le s\le n,
$$
as the {\it structure polynomial matrices} of $g$.

\begin{defi}\label{def-sv}\rm
A Virasoro element $L\in \mathfrak{gc}_N$ is {\it standard} if every structure polynomial matrix $P_L^{[s]}(\partial)$ of
$L$ admits a decomposition $P_L^{[s]}(\partial)=f_L^{[s]}(\partial)P_L^{[s]}$, where $f_L^{[s]}(\partial)\in\C[\partial]$ and $P_L^{[s]}\in M_N$.
\end{defi}

Denote by ${\rm Vir}_N$ (resp., ${\rm Vir}_N^{\rm can}$ and ${\rm Vir}_N^{\rm std}$) the set of all (resp., canonical and standard)
Virasoro elements of $\mathfrak{gc}_N$.
By Theorem~\ref{prop-vir-ele-1},
$${\rm Vir}_1^{\rm can}= {\rm Vir}_1^{\rm std} = {\rm Vir}_1.$$
By Definition~\ref{def-sv}, it is clear that all canonical Virasoro elements of $\mathfrak{gc}_N$ are standard,
namely, ${\rm Vir}_N^{\rm can}\subseteq {\rm Vir}_N^{\rm std}$.
Obviously, any element of $\mathfrak{gc}_N$ of degree zero can not be a Virasoro element, let alone a standard one.
For those of degree one, we have the following classification, which in particular implies that ${\rm Vir}_N^{\rm can}\varsubsetneq {\rm Vir}_N^{\rm std}$ for  $N\ge2$.

\begin{theo}\label{prop-vir-ele-N}
An element of $\mathfrak{gc}_N$ of degree one is a standard Virasoro element if and only if it has the form
$(a\partial+b)J^0_{A_0}+J^1_{A_1}$, where $a,b\in\C$, $A_0, A_1\in M_N$, satisfying $A_1\ne \mathbf{0}$, $A_1 A_0=A_0$,
and $(A_1-a A_0)(A_1-I_N)=\mathbf{0}$.
\end{theo}

\ni{\it Proof.}\ \ The sufficiency can be checked directly.

Next, we prove the necessity. Assume that $L\in {\rm Vir}_N^{\rm std}$ has degree one.
By Definition~\ref{def-sv}, we can write
$$
L=f_{0}(\partial)J^0_{A_{0}}+f_{1}(\partial)J^1_{A_{1}},
$$
where $f_{0}(\partial),f_{1}(\partial)\in \C[\partial]$, $A_{0},A_{1}\in M_N$ and $f_{1}(\partial)\ne 0$, $A_{1}\ne \mathbf{0}$.
Since $[L\, {}_\lambda \, L]=(\partial+2\lambda)L$, we have
\begin{eqnarray}
\nonumber \mathbf{L}_1+\mathbf{L}_2+\mathbf{L}_3 &\!\!\!=\!\!\!&  (\partial+2\lambda)f_0(\partial)J^0_{A_0}+(\partial+2\lambda)f_1(\partial)J^1_{A_1}, \quad \mbox{where} \\
\nonumber \mathbf{L}_1 &\!\!\!=\!\!\!& \left(\lambda f_0(-\lambda)f_1(\partial+\lambda)+(\partial+\lambda)f_1(-\lambda)f_0(\partial+\lambda)\right)J^0_{A_1A_0}, \\
\nonumber \mathbf{L}_2 &\!\!\!=\!\!\!& (f_0(-\lambda)f_1(\partial+\lambda)-f_1(-\lambda)f_0(\partial+\lambda))J^1_{[A_0,\,A_1]}, \\
\label{equ-LL-new} \mathbf{L}_3 &\!\!\!=\!\!\!& (\partial+2\lambda) f_1(-\lambda)f_1(\partial+\lambda)J^1_{A^2_1}.
\end{eqnarray}
Here, $[A_0, A_1]=A_0 A_1-A_1 A_0$.

{\bf Case 1.} Assume $f_0(\partial)=0$.
By \eqref{equ-LL-new}, we see that $f_1(-\lambda)f_1(\partial+\lambda)=f_1(\partial)$ and $A^2_1=A_1$.
Hence, $f_1(\partial)=1$, and thus $L=J^1_{A_{1}}$, which is of the required form (with $A_0=\mathbf{0}$).

{\bf Case 2.} Assume $f_0(\partial)\ne 0$. Comparing the terms $J^0_{*}$'s in \eqref{equ-LL-new}, we obtain
\begin{eqnarray}
\label{equ-case2-1} \lambda f_0(-\lambda)f_1(\partial+\lambda)+(\partial+\lambda)f_1(-\lambda)f_0(\partial+\lambda)&\!\!\!=\!\!\!& (\partial+2\lambda)f_0(\partial),\\
\label{equ-case2-2}  A_1 A_0 &\!\!\!=\!\!\!& A_0.
\end{eqnarray}

{\bf Subcase 2.1.} Assume $[A_0, A_1] = \mathbf{0}$.
Comparing the terms $J^1_{*}$'s in \eqref{equ-LL-new}, as in Case~1, we obtain $f_1(\partial)=1$ and $A^2_1=A_1$.
Then, using the same arguments as in the proof of Theorem~\ref{prop-vir-ele-1}, by \eqref{equ-case2-1} we obtain
$f_0(\partial)=a\partial+b$, where $a,b\in\C$. 
Finally, by \eqref{equ-case2-2} and relations $[A_0, A_1] = \mathbf{0}$ and $A^2_1=A_1$, one can easily check that $(A_1-a A_0)(A_1-I_N)=\mathbf{0}$, as required.

{\bf Subcase 2.2.} Assume $f_0(-\lambda)f_1(\partial+\lambda)=f_1(-\lambda)f_0(\partial+\lambda)$.
Comparing the terms $J^1_{*}$'s in \eqref{equ-LL-new}, we again obtain $f_1(\partial)=1$ and $A^2_1=A_1$.
Then $f_0(-\lambda)=f_0(\partial+\lambda)$, which implies $f_0(\partial)=b\in\C$.
By \eqref{equ-case2-2} and the relation $A^2_1=A_1$, we see that $L$ has the required form (with $a=0$).

{\bf Subcase 2.3.} Assume $[A_0, A_1] \ne \mathbf{0}$ and $f_0(-\lambda)f_1(\partial+\lambda) \ne f_1(-\lambda)f_0(\partial+\lambda)$.
Comparing the terms $J^1_{*}$'s in \eqref{equ-LL-new}, we obtain
$$
(f_0(-\lambda)f_1(\partial+\lambda)-f_1(-\lambda)f_0(\partial+\lambda))J^1_{[A_0, A_1]}
+(\partial+2\lambda) f_1(-\lambda)f_1(\partial+\lambda)J^1_{A^2_1} = (\partial+2\lambda)f_1(\partial)J^1_{A_1}.
$$
Then there are two possibilities. One possibility is the system
\begin{eqnarray}
\label{subsub1-1} (\partial+2\lambda) f_1(-\lambda)f_1(\partial+\lambda) &\!\!\!=\!\!\!& k_1 (f_0(-\lambda)f_1(\partial+\lambda)-f_1(-\lambda)f_0(\partial+\lambda)), \\
\label{subsub1-2} f_0(-\lambda)f_1(\partial+\lambda)-f_1(-\lambda)f_0(\partial+\lambda) &\!\!\!=\!\!\!& k_2 (\partial+2\lambda)f_1(\partial), \\
\label{subsub1-3} A_1 &\!\!\!=\!\!\!& k_2 ([A_0, A_1]+k_1 A_1^2),  
\end{eqnarray}
where $k_1, k_2\in\C$ and $k_1, k_2\ne 0$.
By \eqref{subsub1-1} and \eqref{subsub1-2}, we see that $f_1(-\lambda)f_1(\partial+\lambda)=k_1 k_2 f_1(\partial)$, which implies $f_1(\partial)=k_1 k_2$.
Substituting this back into \eqref{subsub1-1} (or \eqref{subsub1-2}), and then taking $\lambda=0$, we obtain $f_0(\partial)=-k_2 \partial +b$, where $b=f_0(0)\in\C$.
Then, comparing the coefficients of $\partial^2$ in \eqref{equ-case2-1}, we see that $k_1 k_2=1$, and so $f_1(\partial)=1$.
Finally, by \eqref{equ-case2-2}, \eqref{subsub1-3} and $k_1 k_2=1$, one can check that $(A_1+ k_2 A_0)(A_1-I_N)=\mathbf{0}$, as required.
The other possibility is the system
\begin{eqnarray}
\label{subsub2-1} A_1^2 &\!\!\!=\!\!\!& k_1 [A_0, A_1], \\
\label{subsub2-2} [A_0, A_1] &\!\!\!=\!\!\!& k_2 A_1, \\
\label{subsub2-3} (\partial+2\lambda)f_1(\partial) &\!\!\!=\!\!\!& k_2 (f_0(-\lambda)f_1(\partial\!+\!\lambda)\!-\!f_1(-\lambda)f_0(\partial\!+\!\lambda)\!+\!k_1 (\partial\!+\!2\lambda) f_1(-\lambda)f_1(\partial\!+\!\lambda)),  
\end{eqnarray}
where $k_1, k_2\in\C$ and $k_2\ne 0$.
Taking $\lambda=0$ in \eqref{equ-case2-1}, we have $\partial f_1(0)f_0(\partial)=\partial f_0(\partial)$, which gives $f_1(0)=1$.
Then, taking $\lambda=0$ in \eqref{subsub2-3}, we obtain
\begin{equation}\label{f0-form}
f_0(\partial)=(a\partial+b)f_1(\partial), \quad \mbox{where}\quad a=k_1-\frac{1}{k_2},\ b=f_0(0).
\end{equation}
Substituting this into \eqref{equ-case2-1}, we obtain 
\begin{equation}\label{f1-form}
(\partial+2\lambda)(a\partial+b)(f_1(-\lambda)f_1(\partial+\lambda)- f_1(\partial))=0.
\end{equation}
Since $f_0(\partial)\ne 0$, by \eqref{f0-form}, we must have $a\partial+b\ne 0$. Then, by \eqref{f1-form}, we have $f_1(-\lambda)f_1(\partial+\lambda)=f_1(\partial)$,
which gives $f_1(\partial)=1$. By \eqref{f0-form}, $f_0(\partial)=a\partial+b$.
Finally, by \eqref{equ-case2-2}, \eqref{subsub2-1}, \eqref{subsub2-2} and $a=k_1-\frac{1}{k_2}$, one can check that $(A_1-a A_0)(A_1-I_N)=\mathbf{0}$, as required.
\QED

\vskip10pt

Taking $a=0$ in Theorem~\ref{prop-vir-ele-N}, we obtain a Virasoro element $L=b J^0_{A_0}+J^1_{A_1}$, where $A_1 A_0=A_0$,
and $A_1^2=A_1\ne \mathbf{0}$. The conditions $A_1 A_0=A_0$ and $A^2_1=A_1$ imply that ${\rm Im}\,A_0\subseteq {\rm Ker}\,(A_1-I_N)={\rm Im}\,A_1$, 
and thus $A_0=A_1 B$ for some $B\in M_N$. Hence, 
$$L=b J^0_{A_1 B}+J^1_{A_1}.$$
Motivated by the above form, we can construct standard Virasoro elements of $\mathfrak{gc}_N$ ($N\ge 2$) of higher degree;
details of verification are omitted.

\begin{prop}\label{std-vir-ele-higher}
Let $N, k\ge 2$, $a_i\in\C$ and $A,B_i\in M_N$, where $2\le i\le k$.
If $A^2=A\ne \mathbf{0}$ and $AB_i A=\mathbf{0}$, then
$$J^1_A+\sum_{i=2}^k a_i J^i_{AB_i}\in\mathfrak{gc}_N$$
is a standard Virasoro element of $\mathfrak{gc}_N$ (note that all structure matrices are constant).
This, in particular, implies that there exist standard Virasoro elements of $\mathfrak{gc}_N$ of any positive degree.
\end{prop}

\subsection{Non-standard Virasoro elements of $\mathfrak{gc}_N$}

Letting $A_0 A_1=A_0$ in Theorem~\ref{prop-vir-ele-N}, we see that $A^2_1=A_1$. As above, the conditions $A_1 A_0=A_0$ and $A^2_1=A_1$ give that $A_0=A_1 B$ for some $B\in M_N$.
Furthermore, the condition $A_0 A_1 = A_0$ implies that $A_0=A_1 B A_1$. Namely, the resulting Virasoro element is
\begin{equation}\label{Vir-subcase2.1}
L=(a\partial + b) J^0_{A_1 B A_1}+J^1_{A_1}.
\end{equation}
Based on this form and Proposition~\ref{std-vir-ele-higher}, we can further construct
non-standard Virasoro elements of $\mathfrak{gc}_N$ ($N\ge 2$).
Hence, ${\rm Vir}_N^{\rm std}\varsubsetneq {\rm Vir}_N$ for  $N\ge2$.

\begin{prop}\label{non-std-vir-ele}
Let $N,k,\ell\ge 2$, $a_i,b_j\in\C$ and $A_i, B_i, C, D_j\in M_N$, where $1\le i\le k$, $2\le j\le \ell$.
If $A_i^2=A_i\ne \mathbf{0}$, $A_i B_i A_i\ne \mathbf{0}$, and $A_i A_j=A_j A_i=\mathbf{0}$,
$A_i B_i A_i$ and $A_j B_j A_j$ are not proportional and $a_i\ne a_j$ for $i\ne j$, then
\begin{eqnarray*}
&&T_1:=J^1_{\sum_{i=1}^k A_i}+J^0_{A_1 B_1 A_1}+\sum_{i=2}^k(\partial+a_i)J^0_{A_i B_i A_i}\in\mathfrak{gc}_N \ \ \mbox{and} \\
&&T_2:=J^1_{\sum_{i=1}^k A_i}+\sum_{i=1}^k(\partial+a_i)J^0_{A_i B_i A_i}\in\mathfrak{gc}_N
\end{eqnarray*}
are non-standard Virasoro elements of $\mathfrak{gc}_N$ of degree one.
If further $N\ge 3$, $C^2=C\ne \mathbf{0}$ and $A_i C=C A_i=C D_j A_i=\mathbf{0}$, then
\begin{eqnarray*}
&&T_3:=J^1_{C+\sum_{i=1}^k A_i}+J^0_{A_1 B_1 A_1}+\sum_{i=2}^k(\partial+a_i)J^0_{A_i B_i A_i}+\sum_{j=2}^\ell b_j J^j_{C D_j}\in\mathfrak{gc}_N \ \ \mbox{and}\\
&&T_4:=J^1_{C+\sum_{i=1}^k A_i}+\sum_{i=1}^k(\partial+a_i)J^0_{A_i B_i A_i}+ \sum_{j=2}^\ell b_j J^j_{C D_j}\in\mathfrak{gc}_N
\end{eqnarray*}
are non-standard Virasoro elements of $\mathfrak{gc}_N$.
This, in particular, implies that there exist non-standard Virasoro elements of $\mathfrak{gc}_N$ of any positive degree for $N\ge 3$.
\end{prop}

\ni{\it Proof.}\ \
The conditions $A_i^2=A_i\ne \mathbf{0}$ and $A_i A_j=A_j A_i=\mathbf{0}$ for $i\ne j$ imply that $\sum_{i=1}^k A_i\ne \mathbf{0}$.
Hence, $T_1$ has degree one. Rewrite $T_1$ as
$$
T_1=\sum_{i=1}^k L_i, \ \mbox{where}\ L_1=J^1_{A_1}+J^0_{A_1 B_1 A_1}\ \mbox{and}\ L_i=J^1_{A_i}+(\partial+a_i)J^0_{A_i B_i A_i}, 2\le i\le k.
$$
By \eqref{Vir-subcase2.1}, $L_i$'s are (standard) Virasoro elements of $\mathfrak{gc}_N$.
The conditions $A_i A_j=A_j A_i=\mathbf{0}$ for $i\ne j$ give $[L_i\, {}_\lambda \, L_j]=0$.
Hence, $T_1$ is a Virasoro element of $\mathfrak{gc}_N$.
Furthermore, the conditions $A_i B_i A_i\ne \mathbf{0}$,
$A_i B_i A_i$ and $A_j B_j A_j$ are not proportional and $a_i\ne a_j$ for $i\ne j$ imply that $T_1$ is non-standard.
Similarly, $T_2$ is also a non-standard Virasoro element of $\mathfrak{gc}_N$ of degree one.

The conditions $C^2=C\ne \mathbf{0}$ and $A_i C=C A_i=\mathbf{0}$ imply that $C+\sum_{i=1}^k A_i\ne \mathbf{0}$.
Hence, the degree of $T_3$ is at least one. Rewrite $T_3$ as
$$
T_3=T_1+L_3, \ \mbox{where}\ L_3=J^1_{C}+\sum_{j=2}^\ell b_j J^j_{C D_j}.
$$
By Proposition~\ref{std-vir-ele-higher}, $L_3$ is a (standard) Virasoro element of $\mathfrak{gc}_N$.
Using similar arguments as above, we can show that $T_3$ is a non-standard Virasoro element of $\mathfrak{gc}_N$.
Similarly, $T_4$ is also a non-standard Virasoro element of $\mathfrak{gc}_N$.
\QED

\begin{exam}\label{ex-N=3}\rm
Let $N=3$, $k=\ell=2$, $a_1=0$, $a_2=1$, $b_2=1$, and
$$
A_1=B_1=E_{11}+E_{21},\
A_2=B_2=-E_{21}+E_{22},\
C=D_2=E_{33}.
$$
Then
\begin{eqnarray*}
T_3 &\!\!\!=\!\!\!& J^1_{I_3}+J^0_{E_{11}+E_{21}}+(\partial+1) J^0_{-E_{21}+E_{22}}+J^2_{E_{33}},\\
T_4 &\!\!\!=\!\!\!& J^1_{I_3}+\partial J^0_{E_{11}+E_{21}}+(\partial+1) J^0_{-E_{21}+E_{22}}+J^2_{E_{33}}.
\end{eqnarray*}
One can check that they are indeed non-standard Virasoro elements of $\mathfrak{gc}_3$ of degree two.
However, if we require $N=2$, then $T_3$ and $T_4$ do not exist.
\end{exam}

For any Virasoro element $g\in{\rm Vir}_N$, there is an embedding from $\CVir$ into $\mathfrak{gc}_N$:
$$
\theta_g: \ \CVir\rightarrow \mathfrak{gc}_N, \quad L\mapsto g.
$$
Consider the {\it standard module} $V=\C[\partial]^N$ of $\mathfrak{gc}_N$ with $\lambda$-actions (cf.~\eqref{gcN-actions})
$$
J^n_A\,{}_\lambda\, v=(\partial+\lambda)^n Av, \ \ v\in \C^N.
$$
The embedding map $\theta_g$ and the standard module $V$ establish a conformal module $V_g=\C[\partial]^N$ over $\CVir$ with the $\lambda$-action
$L\,{}_\lambda\, v=\theta_g(L)\,{}_\lambda\, v$ for $v\in \C^N$.

\begin{rema}\label{new-vir-modules}\rm
Taking a canonical Virasoro element $g=\Pi(L_{a,b})\in {\rm Vir}_N^{\rm can}$, we obtain
$$V_{\Pi(L_{a,b})}\cong (M^{\CVir}_{1-a,b})^{(N)}.$$
Taking Virasoro elements in ${\rm Vir}_N$ ($N\ge 2$) to be those in Theorem~\ref{prop-vir-ele-N}, Proposition~\ref{std-vir-ele-higher}
and Proposition~\ref{non-std-vir-ele}, we can obtain a huge number of new Virasoro conformal modules.
For example, taking a standard Virasoro element $g=J^1_A+\sum_{i=2}^k a_i J^i_{AB_i}\in {\rm Vir}_N^{\rm std}$ in Proposition~\ref{std-vir-ele-higher},  we obtain
$$
V_g=\C[\partial]^N: \quad L\,{}_\lambda\, v=(\partial+\lambda)Av+\sum_{i=2}^k a_i (\partial+\lambda)^i AB_i v, \ \ v\in \C^N,
$$
where $k\ge2$, $a_i\in\C$, $A,B_i\in M_N$ and $A^2=A\ne \mathbf{0}$, $AB_i A=\mathbf{0}$.
\end{rema}

\section{Semisimplicity criteria for conformal modules over $\mathfrak{gc}_1$}

We now return to find semisimplicity criteria for finite conformal modules over $\mathfrak{gc}_N$.
In this section, we consider the case $N=1$.
Some results on the representation theory of the Heisenberg-Virasoro conformal algebra will be used.

\subsection{Heisenberg-Virasoro conformal subalgebra of $\mathfrak{gc}_1$}

By \eqref{gcN-brackets} with $N=1$, it is clear that
$$
[J^1\, {}_\lambda \, J^1]=(\partial+2\lambda) J^1,\quad [J^1\, {}_\lambda \, J^0]=(\partial+\lambda) J^0,
\quad [J^0\, {}_\lambda \, J^0]=0.
$$
In other words,  $\mathfrak{HV}=\C[\partial]J^0\oplus\C[\partial]J^1$ is a Heisenberg-Virasoro conformal subalgebra of $\mathfrak{gc}_1$.
As a direct corollary of Theorem~\ref{prop-vir-ele-1}, any Virasoro element of $\mathfrak{HV}$ has the form $(a\partial+b)J^0+J^1$, where $a,b\in\C$.
We have the following classification results (see, e.g., \cite{SXY2018}) on representations of $\mathfrak{HV}$.

\begin{lemm}\label{HV-rep}
Any non-trivial free conformal module of rank one over $\mathfrak{HV}$ has the form
$M^{\mathfrak{HV}}_{\D,\a,\beta}=\C[\partial]v$ with $\lambda$-actions
\begin{equation*}
J^1\,{}_\lambda\, v=(\partial+\D\lambda+\a)v, \quad J^0\,{}_\lambda\, v=\beta v,
\end{equation*}
where $\D,\a,\beta\in\C$. 
Furthermore, $M^{\mathfrak{HV}}_{\D,\a,\beta}$ is irreducible if and only if $\D\ne 0$ or $\beta\ne 0$,
and all non-trivial finite irreducible conformal modules over $\mathfrak{HV}$ are of this kind.
\end{lemm}

Using Lemma~\ref{HV-rep}, we can establish a relation between the regularity of the actions of Virasoro elements
in conformal $\mathfrak{gc}_1$-modules with conformal $\mathfrak{HV}$-modules.

\begin{lemm}\label{case1-step1}
Let $V$ be a finite conformal module over $\mathfrak{gc}_1$.
Then there exists a pair of different Virasoro elements, whose actions on a $\C[\partial]$-basis of $V$ are regular
if and only if $V$ as a $\mathfrak{HV}$-module is a direct sum of conformal modules of rank one.
\end{lemm}

\ni{\it Proof.}\ \
The sufficiency follows easily from Theorem~\ref{prop-vir-ele-1} and Lemma~\ref{HV-rep}.

Next we prove the necessity. Assume that $V$ has rank $k$.
Recall by Theorem~\ref{prop-vir-ele-1} that any Virasoro element of $\mathfrak{gc}_1$ has the form $L_{a,b}=(a\partial+b)J^0+J^1$, where $a,b\in\C$.
Suppose that $L_{a_1,b_1}$ and $L_{a_2,b_2}$ are two different Virasoro elements, whose actions on a $\C[\partial]$-basis
$Y=\{v_i\,|\,1\le i\le k\}$ are regular. By Definition~\ref{def-ra}, we have
\begin{eqnarray}
\label{vir-act-1} L_{a_1,b_1}\,{}_\lambda\, v_i &\!\!\!=\!\!\!& \left(\partial+\D^{(i)}_{a_1,b_1}\lambda+\a^{(i)}_{a_1,b_1}\right)v_i,\\
\label{vir-act-2} L_{a_2,b_2}\,{}_\lambda\, v_i &\!\!\!=\!\!\!& \left(\partial+\D^{(i)}_{a_2,b_2}\lambda+\a^{(i)}_{a_2,b_2}\right)v_i,
\end{eqnarray}
where $\D^{(i)}_{a_1,b_1},\D^{(i)}_{a_2,b_2},\a^{(i)}_{a_1,b_1},\a^{(i)}_{a_2,b_2}\in\C$. 
Subtracting \eqref{vir-act-2} from \eqref{vir-act-1}, we obtain
\begin{equation}\label{vir-act-1-case1}
\left(-(a_1-a_2)\lambda+(b_1-b_2)\right)J^0\,{}_\lambda\, v_i=\left(\left(\D^{(i)}_{a_1,b_1}-\D^{(i)}_{a_2,b_2}\right)\lambda+\left(\a^{(i)}_{a_1,b_1}-\a^{(i)}_{a_2,b_2}\right)\right)v_i.
\end{equation}
If $a_1\ne a_2$, then by \eqref{vir-act-1-case1} we have
$$
-(a_1-a_2)\lambda+(b_1-b_2)\mid \left(\D^{(i)}_{a_1,b_1}-\D^{(i)}_{a_2,b_2}\right)\lambda+\left(\a^{(i)}_{a_1,b_1}-\a^{(i)}_{a_2,b_2}\right),
$$
and thus there exists some $\beta_i\in\C$ such that $J^0\,{}_\lambda\, v_i=\beta_i v_i\in \C[\lambda]v_i$.
If $a_1=a_2$, then $b_1\ne b_2$ (since $L_{a,b_1}\ne L_{a,b_2}$).
By \eqref{vir-act-1-case1}, we still have $J^0\,{}_\lambda\, v_i\in \C[\lambda]v_i$.
Then, by \eqref{vir-act-1}, we have $J^1\,{}_\lambda\, v_i\in \C[\partial,\lambda]v_i$.
Hence, $\C[\partial]v_i$ is a conformal module of rank one over $\mathfrak{HV}$, and thus $V$ is a direct sum of these conformal modules.
\QED

\subsection{Proof of part (1) of the Main Theorem}

Now we consider the proof of part (1) of the Main Theorem.
The necessity follows easily from the structure of finite irreducible conformal modules over $\mathfrak{gc}_1$ (see Example~2)
and Theorem~\ref{prop-vir-ele-1}.
In this subsection, we mainly prove the sufficiency:

\begin{theo}\label{thm-1}
Let $V$ be a finite conformal module over $\mathfrak{gc}_1$.
If there exists a pair of different Virasoro elements, whose actions on a $\C[\partial]$-basis of $V$ are regular, then $V$ is semisimple.
\end{theo}

\ni{\it Proof.}\ \
Assume that $V$ has rank $k$. 
By Lemmas~\ref{HV-rep} and \ref{case1-step1}, there exists a $\C[\partial]$-basis $Y=$\linebreak$\{v_i\,|\,1\le i\le k\}$ of $V$
such that
\begin{equation}\label{J0J1}
J^0\,{}_\lambda\, v_i=\beta_i v_i, \quad J^1\,{}_\lambda\, v_i=(\partial+\D_i\lambda+\alpha_i) v_i,
\end{equation}
where $\D_i, \alpha_i,\beta_i\in\C$.
For $n\in\Z_+$ and $1\le i\le k$, assume that
\begin{equation*}
J^n\,{}_\lambda\, v_i=\sum_{j=1}^k f_{ij}^{[n]}(\partial,\lambda)v_j,
\end{equation*}
where $f_{ij}^{[n]}(\partial,\lambda)\in\C[\partial,\lambda]$.
For $m,n\in\Z_+$, introduce the {\it conformal commutator}
$$
T^{m,n}(\lambda,\mu):=J^m_\lambda J^n_\mu-J^n_\mu J^m_\lambda=[J^m\, {}_\lambda \, J^n]_{\lambda+\mu},
$$
where the $\lambda$-bracket is given by \eqref{gcN-brackets} with $N=1$.
Next we shall frequently use the above two expressions of conformal commutators.
First, applying the operator $T^{0,2}(\lambda,\mu)$ on $v_i$, we obtain
$$
\sum_{j=1}^k \beta_j f_{ij}^{[2]}(\partial+\lambda,\mu)v_j- \sum_{j=1}^k \beta_i f_{ij}^{[2]}(\partial,\mu)v_j=
2\lambda(\partial+\D_i(\lambda+\mu)+\alpha_i)v_i-\beta_i\lambda^2 v_i.
$$
Comparing the coefficients of $v_i$, we have
\begin{equation}\label{J0J2-vi}
\beta_i\left(f_{ii}^{[2]}(\partial+\lambda,\mu)-f_{ii}^{[2]}(\partial,\mu)\right)=2\lambda(\partial+\D_i(\lambda+\mu)+\alpha_i)-\beta_i\lambda^2.
\end{equation}
Note first that $\beta_i\ne0$ (otherwise \eqref{J0J2-vi} gives a contradiction).
Rewrite \eqref{J0J2-vi} as
$$
\frac{f_{ii}^{[2]}(\partial+\lambda,\mu)-f_{ii}^{[2]}(\partial,\mu)}{\lambda}=\frac{2}{\beta_i}(\partial+\D_i(\lambda+\mu)+\alpha_i)-\lambda.
$$
Taking $\lambda\rightarrow 0$, we obtain $\frac{d}{d\partial}f_{ii}^{[2]}(\partial,\mu)=\frac{2}{\beta_i}(\partial+\D_i \mu +\alpha_i)$.
Hence,
$$
f_{ii}^{[2]}(\partial,\lambda)=\frac{1}{\beta_i}\partial^2+\frac{2}{\beta_i}\partial(\D_i\lambda+\alpha_i)+\phi_{i}^{[2]}(\lambda),
$$
where $\phi_{i}^{[2]}(\lambda)\in\C[\lambda]$.
Furthermore, applying the operator $T^{0,3}(\lambda,\mu)$ on $v_i$, and then comparing the coefficients of $v_i$, we obtain
\begin{equation}\label{J0J3-vi}
\beta_i\left(f_{ii}^{[3]}(\partial+\lambda,\mu)-f_{ii}^{[3]}(\partial,\mu)\right)=3\lambda f_{ii}^{[2]}(\partial,\lambda+\mu)-3\lambda^2(\partial+\D_i(\lambda+\mu)+\alpha_i)+\beta_i\lambda^3.
\end{equation}
By standard  calculus techniques, we obtain the preliminary form of $f_{ii}^{[3]}(\partial,\lambda)$:
$$
f_{ii}^{[3]}(\partial,\lambda)=\frac{1}{\beta^2_i}\partial^3+\frac{3}{\beta^2_i}\partial^2(\D_i\lambda+\alpha_i)+\frac{3}{\beta_i}\partial\phi_{i}^{[2]}(\lambda)+\phi_{i}^{[3]}(\lambda),
$$
where $\phi_{i}^{[3]}(\lambda)\in\C[\lambda]$.
Substituting this back into \eqref{J0J3-vi}, we obtain
\begin{eqnarray}
\label{phi-2} 3\beta_i\left(\phi_{i}^{[2]}(\lambda+\mu)-\phi_{i}^{[2]}(\mu)\right) &\!\!\!=\!\!\!& (1+3\D_i\beta_i-\beta_i^2)\lambda^2+3(1+\beta_i)(\D_i\mu+\alpha_i)\lambda, \\
\label{relation-1} \beta_i &\!\!\!=\!\!\!& 2\D_i-1.
\end{eqnarray}
Applying standard calculus techniques further on \eqref{phi-2}, we obtain
$$
\phi_{i}^{[2]}(\lambda)=\left(1+\frac{1}{\beta_i}\right)\left(\frac{1}{2}\D_i\lambda^2+\alpha_i\lambda\right)+c_i,
$$
where $c_i\in\C$. Substituting this back into \eqref{phi-2}, and then comparing the coefficients of $\lambda^2$, we obtain another relation between $\beta_i$ and $\D_i$:
\begin{equation}\label{relation-2}
2\beta_i^2-3\D_i\beta_i+3\D_i-2=0.
\end{equation}
The system \eqref{relation-1} and \eqref{relation-2} imply that there exists a partition of the set $K=\{i\in\Z\,|\,1\le i\le k\}$:
$$
K=\{i\in\Z\,|\,1\le i\le k\}=K_1\bigcup K_2,
$$
such that
\begin{equation}\label{partition}
\left\{\begin{array}{lll}
\D_i=1,\ \beta_i=1, & \mbox{\ if\ \ } i\in K_1,\\[3pt]
\D_i=0,\ \beta_i=-1, & \mbox{\ if\ \ } i\in K_2.
\end{array}\right.
\end{equation}

Next, we prove 

{\bf Claim 1.} For any $n\in\Z_+$, $f_{ij}^{[n]}(\partial,\lambda)=0$ if $i\ne j$.

The cases for $n=0,1$ have been guaranteed by \eqref{J0J1}.
Assume that $n\ge 2$, and the claim holds for $m<n$. Next, we consider the case $n$.
Applying the operator $T^{0,n}(\lambda,\mu)$ on $v_i$ together with the inductive hypothesis, and then comparing the coefficients of $v_j (j\ne i)$,
we obtain
\begin{equation}\label{fijn-claim1}
f_{ij}^{[n]}(\partial+\lambda,\mu)=\frac{\beta_i}{\beta_j}f_{ij}^{[n]}(\partial,\mu).
\end{equation}

If $i\in K_1$, $j\in K_2$ (or $i\in K_2$, $j\in K_1$), then $\beta_i=-\beta_j$ by \eqref{partition}.
Taking $\lambda=0$ in \eqref{fijn-claim1}, we see that Claim~1 holds.

If $i,j\in K_1$ (or $i,j\in K_2$), then $\D_i=\D_j$ and $\beta_i=\beta_j$  by \eqref{partition}.
Then \eqref{fijn-claim1} gives $f_{ij}^{[n]}(\partial,\lambda)=f_{ij}^{[n]}(0,\lambda)$.
Applying the operator $T^{1,n}(\lambda,\mu)$ on $v_i$, and then comparing the coefficients of $v_j (j\ne i)$, we obtain
\begin{equation*}
(\alpha_j-\alpha_i-\mu)f_{ij}^{[n]}(0,\mu)=(n\lambda-\mu)f_{ij}^{[n]}(0,\lambda+\mu),
\end{equation*}
which implies that $f_{ij}^{[n]}(0,\lambda)=0$, and thus Claim~1 holds.

Furthermore, we prove 

{\bf Claim 2.} For $n\in\Z_+$, we have
\begin{equation*}\label{fii}
f_{ii}^{[n]}(\partial,\lambda)=
\left\{\begin{array}{lll}
(\partial+\lambda+\alpha_i)^n, & \mbox{\ if\ \ } i\in K_1,\\[3pt]
-(-\partial-\alpha_i)^n, & \mbox{\ if\ \ } i\in K_2.
\end{array}\right.
\end{equation*}

We only prove the case for $i\in K_1$; the case for $i\in K_2$ can be proven similarly.
The cases for $n=0,1$ have been given by \eqref{J0J1} and \eqref{partition}.
Assume that $n\ge 2$ and the claim holds for $m<n$. Next, we consider the case $n$.
Applying the operator $T^{0,n}(\lambda,\mu)$ on $v_i$ together with the inductive hypothesis, and then comparing the coefficients of $v_i$,
 we can derive that
$$
f_{ii}^{[n]}(\partial+\lambda,\mu)-f_{ii}^{[n]}(\partial,\mu)=(\partial+\lambda+\mu+\alpha_i)^n-(\partial+\mu+\alpha_i)^n.
$$
By standard calculus techniques, we have
\begin{equation}\label{fiin-claim2}
f_{ii}^{[n]}(\partial,\lambda)=(\partial+\lambda+\alpha_i)^n+\phi_{i}^{[n]}(\lambda),
\end{equation}
where $\phi_{i}^{[n]}(\lambda)\in\C[\lambda]$.
Furthermore, applying the operator $T^{1,n}(\lambda,\mu)$ on $v_i$, we obtain
$\mu \phi_{i}^{[n]}(\mu)=(\mu-n\lambda)\phi_{i}^{[n]}(\lambda+\mu)$, which implies that $\phi_{i}^{[n]}(\lambda)=0$.
Then Claim~2 follows from  \eqref{fiin-claim2}.

By Claims~1 and 2, we have in fact obtained the direct sum decomposition
$$
V=\Bigg(\bigoplus_{i\in K_1}M^{\mathfrak{gc}_1}_{\a_i}\Bigg)\bigoplus\Bigg(\bigoplus_{i\in K_2}(M^{\mathfrak{gc}_1}_{-\a_i})^*\Bigg).
$$
Namely, $V$ is semisimple.
\QED

\section{Semisimplicity criteria for conformal modules over $\mathfrak{gc}_N$}

In this section, we find semisimplicity criteria for finite conformal modules over $\mathfrak{gc}_N$ for $N\ge 2$.
Although we have no explicit expression of all (even the standard) Virasoro elements of $\mathfrak{gc}_N$, we shall see that
it suffices to use the canonical Virasoro elements.
The rigidities of homomorphisms and anti-homomorphisms between matrix algebras play a crucial role.
Some combinatorial formulas will be also employed.

\subsection{Rigidities of (anti-)homomorphisms between matrix algebras}

The classical Skolem-Noether theorem (see, e.g., \cite{J1989}) is a fundamental result
in the theory of central simple algebras and representation theory:

\begin{theo}[Skolem-Noether]\label{Skolem-Noether}
Let $\mathcal{A}$ be a simple subalgebra of a finite-dimensional central simple algebra $\mathcal{B}$.
Then any algebra homomorphism of $\mathcal{A}$ into $\mathcal{B}$ can be extended to an inner automorphism of $\mathcal{B}$.
\end{theo}

Using the Skolem-Noether theorem, we can further prove 
the following result on the rigidities of homomorphisms and anti-homomorphisms between matrix algebras.

\begin{lemm}\label{rigidity}
Let $\Phi$ be a non-trivial algebra homomorphism or anti-homomorphism from $M_{n_1}$ to $M_{n_2}$.
Then there exist a positive integer $m\in\Z_{\ge1}$ and an invertible matrix $P\in M_{n_2}$ such that $n_2=m n_1$ and, for $A\in M_{n_1}$,
\begin{equation*}
\Phi(A)=
\left\{\begin{array}{lll}
P(I_m\otimes A)P^{-1}, & \mbox{\ if\ \ } \Phi \ \mbox{is an algebra homomorphism},\\[5pt]
P(I_m\otimes A^T)P^{-1},  & \mbox{\ if\ \ } \Phi \ \mbox{is an algebra anti-homomorphism},
\end{array}\right.
\end{equation*}
where $\otimes$ denotes the Kronecker product of matrices.
\end{lemm}

\ni{\it Proof.}\ \ First, consider the case that $\Phi$ is a non-trivial algebra homomorphism from $M_{n_1}$ to $M_{n_2}$.
Since the kernel ${\rm Ker}\,\Phi$ of $\Phi$ is an ideal of the simple algebra $M_{n_1}$, we must have ${\rm Ker}\,\Phi=\{0\}$, and so $\Phi$ is injective.
The homomorphism $\Phi$ allows us to view $\C^{n_2}$ as an $M_{n_1}$-module with the action $A\cdot v=\Phi(A)v$ for $A\in M_{n_1}$ and $v\in\C^{n_2}$.
By the semisimplicity of representations of $M_{n_1}$, there exists $m\in\Z_{\ge1}$ such that $\C^{n_2}\cong(\C^{n_1})^{(m)}$ (the direct sum of $m$ copies of $\C^{n_1}$).
Hence, $n_2=m n_1$, and so $M_{n_2}\cong M_m\otimes M_{n_1}$, where $\otimes$ denotes the Kronecker product of matrices. Denote by $S$ the standard embedding from $M_{n_1}$ to $M_{n_2}$, i.e., 
$$
S: M_{n_1}\rightarrow M_{n_2}, \quad A\mapsto I_m\otimes A. 
$$
By the Skolem-Noether theorem (Theorem~\ref{Skolem-Noether}), there exists an invertible matrix $P\in M_{n_2}$ such that 
$\Phi(A)=P S(A) P^{-1}=P(I_m\otimes A)P^{-1}$.

Next, consider the case that $\Phi$ is a non-trivial algebra anti-homomorphism from $M_{n_1}$ to $M_{n_2}$.
Denote by $\Theta$ the transpose map: $M_{n_1}\rightarrow M_{n_1}, A\mapsto A^T$.
Then $\Psi:=\Phi\circ\Theta$ is an algebra homomorphism from $M_{n_1}$ to $M_{n_2}$. Applying the above arguments, there exist $m\in\Z_{\ge1}$ and an invertible matrix $P\in M_{n_2}$
such that $\Psi(A)=P(I_m\otimes A)P^{-1}$. Since $\Phi(A)=(\Phi\circ\Theta)(A^T)=\Psi(A^T)$, we have $\Phi(A)=P(I_m\otimes A^T)P^{-1}$.
\QED

\vskip10pt

The following combinatorial formulas can be checked straightforward.

\begin{lemm}\label{formulas}
The following formulas hold:
\begin{eqnarray*}
{\rm(1)} && \sum_{i=0}^{n}\binom{n}{i}(-y)^i(x+y)^{n-i}=x^n;\\
{\rm(2)} && \sum_{i=2}^{n+1}\binom{n+1}{i}(-y)^i(x+y)^{n+1-i} = (n+1)y(x+y)^n-(x+y)^{n+1}+x^{n+1};\\
{\rm(3)} && \sum_{i=2}^{n+1}\binom{n+1}{i}(-y)^i x^{n+1-i} = (n+1)yx^n+(x-y)^{n+1}-x^{n+1}.
\end{eqnarray*}
\end{lemm}

\subsection{Proof of part (2) of the Main Theorem}

Similar to the case $N=1$, the necessity of part (2) of the Main Theorem
follows easily from the structure of finite irreducible conformal modules over $\mathfrak{gc}_N$ (see Example~2), Theorem~\ref{prop-vir-ele-1}
and the definition of canonical Virasoro elements of $\mathfrak{gc}_N$ introduced in Subsection~\ref{sec-cs}.
In this subsection, we mainly prove the sufficiency:

\begin{theo}\label{thm-2}
Let $V$ be a finite conformal module over $\mathfrak{gc}_N$.
If there exists a pair of different canonical Virasoro elements, whose actions on a $\C[\partial]$-basis of $V$ are regular, then $V$ is semisimple.
\end{theo}

\ni{\it Proof.}\ \
Assume that $V$ has rank $k$.
Recall that $\mathfrak{gc}_1$ can be embedded into $\mathfrak{gc}_N$ via the canonical embedding map $\Pi$.
Hence $V$ can be viewed as a conformal module over $\Pi(\mathfrak{gc}_1)\cong \mathfrak{gc}_1$.
By the assumption of this theorem and Theorem~\ref{thm-1}, there exist a $\C[\partial]$-basis $Y=\{v_i\,|\,1\le i\le k\}$ of $V$
and a partition of the index set $K=\{i\in\Z\,|\,1\le i\le k\}$:
$$
K=\{i\in\Z\,|\,1\le i\le k\}=K^{(1)}\bigcup K^{(2)},
$$
such that
\begin{equation*}
J^n_{I_N}\,{}_\lambda\, v_i=
\left\{\begin{array}{lll}
\left(\partial+\lambda+\alpha_i^{(1)}\right)^n v_i, & \mbox{\ if\ \ } i\in K^{(1)},\\[5pt]
-\left(-\partial-\alpha_i^{(2)}\right)^n v_i,  & \mbox{\ if\ \ } i\in K^{(2)},
\end{array}\right.
\end{equation*}
where $\alpha_i^{(1)},\alpha_i^{(2)}\in\C$.
For later use, for $\ell=1,2$, we denote by $\{\beta_1^{(\ell)},\ldots,\beta_{s_\ell}^{(\ell)}\}$ the set of all distinct values of the set
$\{\alpha_i^{(\ell)}\mid i\in K^{(\ell)}\}$, and define
\begin{equation}\label{Ksl}
K_s^{(\ell)}=\{i\in K^{(\ell)}\mid \alpha_i^{(\ell)}=\beta_s^{(\ell)}\}, \quad 1\le s\le s_\ell.
\end{equation}
For $n\in\Z_+$, $A\in M_N$ and $1\le i\le k$, assume that
\begin{equation*}
J^n_A\,{}_\lambda\, v_i=\sum_{j\in K} f_{A;ij}^{[n]}(\partial,\lambda)v_j,
\end{equation*}
where $f_{A;ij}^{[n]}(\partial,\lambda)\in\C[\partial,\lambda]$.
For $m,n\in\Z_+$ and $A,B\in M_N$, introduce the more general {\it conformal commutator}
$$
T^{m,n}_{A,B}(\lambda,\mu):=J^m_A{\;}_\lambda J^n_B{\;}_\mu-J^n_B{\;}_\mu J^m_A{\;}_\lambda=[J^m_A\, {}_\lambda \, J^n_B]{\;}_{\lambda+\mu},
$$
where the $\lambda$-bracket is given by \eqref{gcN-brackets}.
We first claim 

{\bf Claim 1.} There exist further partitions of $K^{(1)}$ and $K^{(2)}$:
\begin{equation}\label{fur-partition}
K^{(1)}=\bigcup_{s=1}^{s_1}K_s^{(1)}, \quad K^{(2)}=\bigcup_{s=1}^{s_2}K_s^{(2)},
\end{equation}
such that
\begin{equation*}
f_{A;ij}^{[n]}(\partial,\lambda)=
\left\{\begin{array}{lll}
c_{ij}^A\left(\partial+\lambda+\beta_s^{(1)}\right)^n, & \mbox{\ if\ \ } i,j\in K_s^{(1)},\ 1\le s\le s_1, \\[5pt]
-c_{ij}^A\left(-\partial-\beta_s^{(2)}\right)^n, & \mbox{\ if\ \ } i,j\in K_s^{(2)},\ 1\le s\le s_2, \\[5pt]
0,  & \mbox{\ otherwise},
\end{array}\right.
\end{equation*}
where $c_{ij}^A\in\C$, and $K_s^{(\ell)}$, $\ell=1,2$, is given by \eqref{Ksl}.

Applying the operator $T^{0,1}_{A,I_N}(\lambda,\mu)$ on $v_i$, and then comparing the coefficients of $v_j$, we obtain
\begin{align}
\label{01-AI-1} \left(\partial\!+\!\lambda\!+\!\mu\!+\!\alpha_i^{(1)}\right)f_{A;ij}^{[0]}(\partial,\lambda)
\!-\!\left(\partial\!+\!\mu\!+\!\alpha_j^{(1)}\right)f_{A;ij}^{[0]}(\partial+\mu,\lambda) & =  \lambda f_{A;ij}^{[0]}(\partial,\lambda+\mu), \ i,j\in K_1, \\
\label{01-AI-2} \left(\partial+\lambda+\alpha_i^{(2)}\right)f_{A;ij}^{[0]}(\partial,\lambda)
-\left(\partial+\alpha_j^{(2)}\right)f_{A;ij}^{[0]}(\partial+\mu,\lambda) & =  \lambda f_{A;ij}^{[0]}(\partial,\lambda+\mu), \ i,j\in K_2.
\end{align}
In \eqref{01-AI-1}, if $\alpha_i^{(1)}\ne \alpha_j^{(1)}$, by taking $\mu=0$, we obtain $f_{A;ij}^{[0]}(\partial,\lambda)=0$;
if $\alpha_i^{(1)}=\alpha_j^{(1)}$, the solution must have the form $f_{A;ij}^{[0]}(\partial,\lambda)=c_{ij}^A\in\C$ (see, e.g., \cite[Corollary ~4.2]{X2022}).
Hence there exists a partition of $K^{(1)}$ as in \eqref{fur-partition} such that
\begin{equation}\label{f0-K1}
f_{A;ij}^{[0]}(\partial,\lambda)=
\left\{\begin{array}{lll}
c_{ij}^A, & \mbox{\ if\ \ } i,j\in K_s^{(1)},\ 1\le s\le s_1, \\[3pt]
0,  & \mbox{\ if\  \ } i\in K_s^{(1)}, j\in K_t^{(1)},\ 1\le s\ne t\le s_1.
\end{array}\right.
\end{equation}
Similarly, by \eqref{01-AI-2}, there exists a partition of $K^{(2)}$ as in \eqref{fur-partition} such that
\begin{equation}\label{f0-K2}
f_{A;ij}^{[0]}(\partial,\lambda)=
\left\{\begin{array}{lll}
-c_{ij}^A, & \mbox{\ if\ \ } i,j\in K_s^{(2)},\ 1\le s\le s_2, \\[3pt]
0,  & \mbox{\ if\  \ } i\in K_s^{(2)}, j\in K_t^{(2)},\ 1\le s\ne t\le s_2.
\end{array}\right.
\end{equation}
Note here that the technical reason we write $-c_{ij}^A$ but not $c_{ij}^A$ as in \eqref{f0-K1} is to reserve the minus sign, so as to keep consistent 
with the form \eqref{gcN-actions} of the dual conformal module of $\mathfrak{gc}_N$. 
Applying the operator $T^{0,1}_{I_N,A}(\lambda,\mu)$ on $v_i$, and then comparing the coefficients of $v_j$, we obtain
\begin{eqnarray*}
f_{A;ij}^{[1]}(\partial,\mu) + f_{A;ij}^{[1]}(\partial+\lambda,\mu) &\!\!\!=\!\!\!& - \lambda f_{A;ij}^{[0]}(\partial,\lambda+\mu), \ i\in K^{(1)}, j\in K^{(2)},\\
f_{A;ij}^{[1]}(\partial,\mu) + f_{A;ij}^{[1]}(\partial+\lambda,\mu) &\!\!\!=\!\!\!& \lambda f_{A;ij}^{[0]}(\partial,\lambda+\mu), \ i\in K^{(2)}, j\in K^{(1)},
\end{eqnarray*}
from which one can easily derive that
\begin{equation}\label{f0f1-K1-K2}
f_{A;ij}^{[0]}(\partial,\lambda)=f_{A;ij}^{[1]}(\partial,\lambda)=0, \ i\in K^{(1)}, j\in K^{(2)}, \ \mbox{or} \ i\in K^{(2)}, j\in K^{(1)}.
\end{equation}
From \eqref{f0-K1}--\eqref{f0f1-K1-K2}, we see that Claim~1 holds for $n=0$.

Assume that $n\ge 1$, and Claim~1 holds for $m<n$. Next, we consider the case $n$.
Applying the operator $T^{0,n+1}_{A,I_N}(\lambda,\mu)$ on $v_i$ together with the inductive hypothesis, and then comparing the coefficients of $v_j$,
we can derive that
\begin{equation}\label{fn-K1}
f_{A;ij}^{[n]}(\partial,\lambda)=
\left\{\begin{array}{lll}
c_{ij}^A\left(\partial+\lambda+\beta_s^{(1)}\right)^n, & \mbox{\ if\ \ } i,j\in K_s^{(1)},\ 1\le s\le s_1, \\[5pt]
0,  & \mbox{\ if\  \ } i\in K_s^{(1)}, j\in K_t^{(1)},\ 1\le s\ne t\le s_1.
\end{array}\right.
\end{equation}
In fact, if $i\in K_s^{(1)}$, we have
\begin{eqnarray*}
T^{0,n+1}_{A,I_N}(\lambda,\mu) (v_i) &\!\!\!=\!\!\!& J^0_A{\;}_\lambda\; J^{n+1}_{I_N}{}_\mu\; v_i-J^{n+1}_{I_N}{}_\mu\; J^0_A{\;}_\lambda\; v_i \\
 &\!\!\!=\!\!\!& \left(\partial+\lambda+\mu+\beta_s^{(1)}\right)^{n+1}\sum_{j\in K_s^{(1)}} c_{ij}^A v_j
-\left(\partial+\mu+\beta_s^{(1)}\right)^{n+1}\sum_{j\in K_s^{(1)}} c_{ij}^A v_j,
\end{eqnarray*}
and
\begin{eqnarray*}
T^{0,n+1}_{A,I_N}(\lambda,\mu) (v_i) &\!\!\!=\!\!\!& (n+1)\lambda J^{n}_{A}{\;}_{\lambda+\mu}\; v_i
-\sum_{m=2}^{n+1}\binom{n+1}{m}(-\lambda)^m J^{n+1-m}_{A}{\;}_{\lambda+\mu}\; v_i \\
 &\!\!\!=\!\!\!& (n+1)\lambda\sum_{j\in K} f_{A;ij}^{[n]}(\partial,\lambda+\mu) v_j \\
 &\!\!\!\!\!\!& -\sum_{j\in K_s^{(1)}}\sum_{m=2}^{n+1}c_{ij}^A\binom{n+1}{m}(-\lambda)^m \left(\partial+\lambda+\mu+\beta_s^{(1)}\right)^{n+1-m} v_j.
\end{eqnarray*}
Comparing the coefficients of $v_j$ with $j\in K_s^{(1)}$, by the formula in Lemma~\ref{formulas}\,(2) with the replacement
$(x,y)\rightsquigarrow(\partial+\mu+\beta_s^{(1)}, \lambda)$, we obtain
$f_{A;ij}^{[n]}(\partial,\lambda)=c_{ij}^A\left(\partial+\lambda+\beta_s^{(1)}\right)^n$.
Comparing the coefficients of $v_j$ with $j\in K_t^{(1)}$, $t\ne s$, we see that $f_{A;ij}^{[n]}(\partial,\lambda)=0$.
Hence, \eqref{fn-K1} holds.
Similarly, applying the formula in Lemma~\ref{formulas}\,(3) with the replacement $(x,y)\rightsquigarrow(-\partial-\beta_s^{(2)}, \lambda)$,
we can derive that
\begin{equation}\label{fn-K2}
f_{A;ij}^{[n]}(\partial,\lambda)=
\left\{\begin{array}{lll}
-c_{ij}^A\left(-\partial-\beta_s^{(2)}\right)^n, & \mbox{\ if\ \ } i,j\in K_s^{(2)},\ 1\le s\le s_2, \\[5pt]
0,  & \mbox{\ if\  \ } i\in K_s^{(2)}, j\in K_t^{(2)},\ 1\le s\ne t\le s_2.
\end{array}\right.
\end{equation}
In addition, applying the operator $T^{0,n+1}_{I_N,A}(\lambda,\mu)$ on $v_i$, and then comparing the coefficients of $v_j$,
we can generalize \eqref{f0f1-K1-K2} to
\begin{equation}\label{fn-K1-K2}
f_{A;ij}^{[n]}(\partial,\lambda)=0, \ i\in K^{(1)}, j\in K^{(2)}, \ \mbox{or} \ i\in K^{(2)}, j\in K^{(1)}.
\end{equation}
Now Claim~1 follows from  \eqref{fn-K1}--\eqref{fn-K1-K2}.

By Claim~1, we have obtained a preliminary direct sum decomposition
\begin{equation}\label{decomposition-pre}
V=\Bigg(\bigoplus_{s=1}^{s_1}V_s^{(1)}\Bigg)\bigoplus\Bigg(\bigoplus_{s=1}^{s_2}V_s^{(2)}\Bigg),
\end{equation}
where $V_s^{(\ell)}$ is a submodule of $V$ with $\C[\partial]$-basis
$$
Y_s^{(\ell)}=\{v_i\,|\,i\in K_s^{(\ell)}\}, \quad 1\le s\le s_{\ell}, \ \ell=1,2.
$$
More precisely, the conformal structure of $V_s^{(\ell)}$ is given by
\begin{equation*}
J^n_A\,{}_\lambda\,Y_s^{(\ell)}=
\left\{\begin{array}{lll}
\left(\partial+\lambda+\beta_s^{(1)}\right)^n Y_s^{(1)} \Phi_s^{(1)}(A), & \mbox{\ if\ \ } \ell=1, \\[5pt]
-\left(-\partial-\beta_s^{(2)}\right)^n Y_s^{(2)} \Phi_s^{(2)}(A), & \mbox{\ if\ \ } \ell=2,
\end{array}\right.
\end{equation*}
where $\Phi_s^{(\ell)}(A)=(c_{ij}^A)^T_{i,j\in K_s^{(\ell)}}\in M_{|K_s^{(\ell)}|}$.
We further claim 

{\bf Claim 2.} We have
\begin{equation*}
\Phi_s^{(\ell)}(AB)=
\left\{\begin{array}{lll}
\Phi_s^{(1)}(A) \Phi_s^{(1)}(B), \ 1\le s\le s_{1}, & \mbox{\ if\ \ } \ell=1, \\[5pt]
\Phi_s^{(2)}(B) \Phi_s^{(2)}(A), \ 1\le s\le s_{2}, & \mbox{\ if\ \ } \ell=2.
\end{array}\right.
\end{equation*}

Applying the operator $T^{m,n}_{A,B}(\lambda,\mu)$ on $Y_s^{(1)}$, and then using the formula in Lemma~\ref{formulas}\,(1), we can derive that
\begin{eqnarray*}
\mathbf{L}^{(1)}_1-\mathbf{L}^{(1)}_2 &\!\!\!=\!\!\!& \mathbf{R}^{(1)}_1-\mathbf{R}^{(1)}_2, \quad \mbox{where} \\
\mathbf{L}^{(1)}_1 &\!\!\!=\!\!\!& \left(\partial+\lambda+\mu+\beta_s^{(1)}\right)^n\left(\partial+\lambda+\beta_s^{(1)}\right)^m Y_s^{(1)} \Phi_s^{(1)}(A) \Phi_s^{(1)}(B), \\
\mathbf{L}^{(1)}_2 &\!\!\!=\!\!\!& \left(\partial+\lambda+\mu+\beta_s^{(1)}\right)^m\left(\partial+\mu+\beta_s^{(1)}\right)^n Y_s^{(1)} \Phi_s^{(1)}(B) \Phi_s^{(1)}(A),\\
\mathbf{R}^{(1)}_1 &\!\!\!=\!\!\!& \left(\partial+\lambda+\mu+\beta_s^{(1)}\right)^n\left(\partial+\lambda+\beta_s^{(1)}\right)^m Y_s^{(1)} \Phi_s^{(1)}(AB), \\
\mathbf{R}^{(1)}_2 &\!\!\!=\!\!\!& \left(\partial+\lambda+\mu+\beta_s^{(1)}\right)^m\left(\partial+\mu+\beta_s^{(1)}\right)^n Y_s^{(1)} \Phi_s^{(1)}(BA),
\end{eqnarray*}
which implies that $\Phi_s^{(1)}(AB)=\Phi_s^{(1)}(A) \Phi_s^{(1)}(B)$. Similarly, applying the operator $T^{m,n}_{A,B}(\lambda,\mu)$ on $Y_s^{(2)}$,
we can derive that
\begin{eqnarray*}
\mathbf{L}^{(2)}_1-\mathbf{L}^{(2)}_2 &\!\!\!=\!\!\!& \mathbf{R}^{(2)}_1-\mathbf{R}^{(2)}_2, \quad \mbox{where} \\
\mathbf{L}^{(2)}_1 &\!\!\!=\!\!\!& \left(-\partial-\lambda-\beta_s^{(2)}\right)^n\left(-\partial-\beta_s^{(2)}\right)^m Y_s^{(2)} \Phi_s^{(2)}(A) \Phi_s^{(2)}(B), \\
\mathbf{L}^{(2)}_2 &\!\!\!=\!\!\!& \left(-\partial-\mu-\beta_s^{(2)}\right)^m\left(-\partial-\beta_s^{(2)}\right)^n Y_s^{(2)} \Phi_s^{(2)}(B) \Phi_s^{(2)}(A),\\
\mathbf{R}^{(2)}_1 &\!\!\!=\!\!\!& \left(-\partial-\lambda-\beta_s^{(2)}\right)^n\left(-\partial-\beta_s^{(2)}\right)^m Y_s^{(2)} \Phi_s^{(2)}(BA), \\
\mathbf{R}^{(2)}_2 &\!\!\!=\!\!\!& \left(-\partial-\mu-\beta_s^{(2)}\right)^m\left(-\partial-\beta_s^{(2)}\right)^n Y_s^{(2)} \Phi_s^{(2)}(AB),
\end{eqnarray*}
which implies that $\Phi_s^{(2)}(AB)=\Phi_s^{(2)}(B) \Phi_s^{(2)}(A)$. Hence, Claim~2 holds.

By Claim~2, $\Phi_s^{(1)}$ is an algebra homomorphism from $M_{N}$ to $M_{|K_s^{(1)}|}$, while
$\Phi_s^{(2)}$ is an algebra anti-homomorphism from $M_{N}$ to $M_{|K_s^{(2)}|}$.
By Lemma~\ref{rigidity}, there exist an integer $m_s^{(\ell)}\in\Z_{\ge1}$ and an invertible matrix $P_s^{(\ell)}\in M_{|K_s^{(\ell)}|}$
such that $|K_s^{(\ell)}|=m_s^{(\ell)} N$ and
\begin{equation*}
\Phi_s^{(\ell)}(A)=
\left\{\begin{array}{lll}
P_s^{(1)}\left(I_{m_s^{(1)}}\otimes A\right)(P_s^{(1)})^{-1}, & \mbox{\ if\ \ } \ell=1, \\[5pt]
P_s^{(2)}\left(I_{m_s^{(2)}}\otimes A^T\right)(P_s^{(2)})^{-1}, & \mbox{\ if\ \ } \ell=2.
\end{array}\right.
\end{equation*}
Under the new $\C[\partial]$-basis $\overline{Y}_s^{(\ell)}=Y_s^{(\ell)}P_s^{(\ell)}$ of $V_s^{(\ell)}$ , we have
\begin{equation*}
J^n_A\,{}_\lambda\,\overline{Y}_s^{(\ell)}=
\left\{\begin{array}{lll}
\left(\partial+\lambda+\beta_s^{(1)}\right)^n \overline{Y}_s^{(1)} \left(I_{m_s^{(1)}}\otimes A\right), & \mbox{\ if\ \ } \ell=1, \\[5pt]
-\left(-\partial-\beta_s^{(2)}\right)^n \overline{Y}_s^{(2)} \left(I_{m_s^{(2)}}\otimes A^T\right), & \mbox{\ if\ \ } \ell=2.
\end{array}\right.
\end{equation*}
Hence, we have the direct sum decompositions
\begin{equation}\label{decomposition-further}
V_s^{(1)}=\left(M^{\mathfrak{gc}_N}_{\beta_s^{(1)}}\right)^{(m_s^{(1)})} \quad \mbox{and}\quad
V_s^{(2)}=\left(\left(M^{\mathfrak{gc}_N}_{-\beta_s^{(2)}}\right)^*\right)^{(m_s^{(2)})},
\end{equation}
which imply that the irreducible conformal module $M^{\mathfrak{gc}_N}_{\beta_s^{(1)}}$ has multiplicity $m_s^{(1)}$ in $V$,
while $\left(M^{\mathfrak{gc}_N}_{-\beta_s^{(2)}}\right)^*$ has multiplicity $m_s^{(2)}$ in $V$.
By \eqref{decomposition-pre} and \eqref{decomposition-further}, $V$ is semisimple.
\QED

\section{Open problems}

Finally, we propose some open problems arising from our study.

{\rm (I)} \ Find semisimplicity criteria for finite conformal modules over infinite Lie conformal subalgebras of $\mathfrak{gc}_N$
with Virasoro elements.

The classical Burnside theorem (see, e.g., \cite{J1989}) states that any subalgebra of $M_N$ that acts irreducibly on $\C^N$ is the whole $M_N$.
Conformal analogue of Burnside theorem in the associative case was formulated and conjectured by Kac \cite{K1999},
and finally confirmed by Kolesnikov \cite{Ko2006}.
In the Lie conformal case, 
it was conjectured in \cite{BKL2003} that any infinite irreducible subalgebra of $\mathfrak{gc}_N$
is conjugate to one of the following subalgebras:
\begin{itemize}\parskip-4pt
\item[\rm(A)] $\mathfrak{gc}_{N,P}$, where $\det P\ne 0$;
\item[\rm(B)] $\mathfrak{oc}_{N,P}$, where $\det P\ne 0$ and $P(-x)=P^T(x)$;
\item[\rm(C)] $\mathfrak{spc}_{N,P}$, where $\det P\ne 0$ and $P(-x)=-P^T(x)$.
\end{itemize}
Around the same time, this conjecture was partially confirmed under certain conditions in \cite{DeK2002,Z2002}.
Among the above list, those with Virasoro elements are more attractive from the viewpoint of physics \cite{DeK2002}.
Irreducible finite conformal modules over these physically important infinite Lie conformal algebras were classified in \cite{BKL2003-2,BM2013}
as corollaries of the classification of finite growth modules.
Our Main Theorem provides semisimplicity criteria for finite conformal modules over $\mathfrak{gc}_N$.
It would be interesting to find the answers for finite conformal modules over these {\it classical} infinite subalgebras of $\mathfrak{gc}_N$.

{\rm (II)} \ Classify standard Virasoro elements of $\mathfrak{gc}_N$ for $N\ge 2$.

Recall from Section~3 that
$$
{\rm Vir}_1^{\rm can}= {\rm Vir}_1^{\rm std} = {\rm Vir}_1 \quad \mbox{and} \quad
{\rm Vir}_N^{\rm can}\varsubsetneq {\rm Vir}_N^{\rm std}\varsubsetneq {\rm Vir}_N \ \ \mbox{for}\ \  N\ge 2.
$$
From Proposition~\ref{non-std-vir-ele}, we see that one can construct non-standard Virasoro elements from standard ones,
and further construct new non-standard Virasoro elements form standard ones and known non-standard ones.
Thus, it is an interesting problem to determine the standard Virasoro elements of $\mathfrak{gc}_N$ for $N\ge 2$.
In Theorem~\ref{prop-vir-ele-N}, we have classified those of degree one.

Besides, note that in the second statement of Proposition~\ref{non-std-vir-ele} there is a restriction $N\ge 3$
(see Example~3 for examples with $N=3$); we have not yet been able to construct non-standard Virasoro elements of $\mathfrak{gc}_2$ of degree $\ge 2$.
It is natural to ask whether all Virasoro elements of $\mathfrak{gc}_2$ of degree $\ge 2$ are standard?

\vskip20pt

\small \ni{\bf Acknowledgement}
We would like to express our sincere gratitude to the referees for their careful reading of the previous version of this manuscript and for their valuable comments and suggestions, 
which have helped us improve the paper significantly.
This work was supported by the National Natural Science Foundation of China (No.~12671036),
and the Fundamental Research Funds for the Central Universities (No.~2024KYJD2006).

\vskip20pt

\large\ni{\bf Declarations}

\vskip15pt

\small\ni{\bf Data availability statement}\ \
Data sharing not applicable to this article as no datasets were generated or analysed during the current study.

\vskip10pt

\ni{\bf Conflict of interest statement}\ \
The author has no competing interests to declare that are relevant to the content of this article.

\end{document}